\documentclass[11pt,leqno,twoside]{article}

\usepackage{amsfonts,amsmath,amsthm,amssymb}
\usepackage{enumerate}
\usepackage{graphics,graphicx,subfigure}

\usepackage{color}
\usepackage{todonotes}
\usepackage{cancel}
\usepackage{url}
\usepackage{hyperref}
\usepackage{makeidx}
\usepackage{showidx}
\usepackage{multicol}        % used for the two-column index
\usepackage{xspace}
\usepackage{stmaryrd}        % for \llbracket, \rrbracket
\usepackage{pifont}          % for \ding{227}
\usepackage{fancybox}        % for oval boxes
\usepackage{bm}
\usepackage{bbm}

 \usepackage{fancyhdr}
\theoremstyle{plain}
 \theoremstyle{definition}
 \newtheorem{lem}{Lemma}
 \newtheorem{defn}[lem]{Definition}
 \newtheorem{thm}[lem]{Theorem}
 \newtheorem{prop}[lem]{Proposition}
 \newtheorem{cor}[lem]{Corollary}
 \newtheorem{notn}[lem]{Notations}
 \newtheorem{pb}[lem]{Problem}
 \newtheorem{form}[lem]{Formulae}
 
 \newtheorem*{rk}{Remark}
 \newtheorem*{com}{Comment}
 \newtheorem*{ex}{Example}
 \theoremstyle{remark}

 \newcommand{\blem}{\begin{lem}}
 \newcommand{\elem}{\end{lem}}
 \newcommand{\bdefn}{\begin{defn}}
 \newcommand{\edefn}{\end{defn}}
 \newcommand{\bthm}{\begin{thm} }
 \newcommand{\ethm}{\end{thm}}
 \newcommand{\bprop}{\begin{prop}}
 \newcommand{\eprop}{\end{prop}}
 \newcommand{\bcor}{\begin{cor}}
 \newcommand{\ecor}{\end{cor}}
 \newcommand{\bnotn}{\begin{notn}}
 \newcommand{\enotn}{\end{notn}}
 \newcommand{\bpb}{\begin{pb}}
 \newcommand{\epb}{\end{pb}}
 \newcommand{\bform}{\begin{form}}
 \newcommand{\eform}{\end{form}}
 \newcommand{\brk}{\begin{rk}}
 \newcommand{\erk}{\end{rk}}
 \newcommand{\bcom}{\begin{com}}
 \newcommand{\ecom}{\end{com}}
 \newcommand{\bex}{\begin{ex}}
 \newcommand{\eex}{\end{ex}}
 \newcommand{\bpf}{\begin{proof}}
 \newcommand{\epf}{\end{proof}}

\newcommand{\ve}{{\bf e}}

\newcommand{\vx}{{\bf x}}
\newcommand{\vy}{{\bf y}}

\newcommand{\cE}{\mathcal{E}}

\newcommand{\cH}{\mathcal{H}}

\newcommand{\cK}{\mathcal{K}}

\newcommand{\cX}{\mathcal{X}}

\newcommand{\bN}{\mathbb{N}}

\newcommand{\bR}{\mathbb{R}}

\newcommand{\be}{\begin{equation}}
\newcommand{\ee}{\end{equation}}
\newcommand{\bal}{\begin{align}}
\newcommand{\eal}{\end{align}}
\newcommand{\ba}{\begin{align*}}
\newcommand{\ea}{\end{align*}}
\newcommand{\bmx}{\begin{matrix}}
\newcommand{\emx}{\end{matrix}}
\newcommand{\bbmx}{\begin{bmatrix}}
\newcommand{\ebmx}{\end{bmatrix}}
\newcommand{\bpmx}{\begin{pmatrix}}
\newcommand{\epmx}{\end{pmatrix}}
\newcommand{\bvmx}{\begin{vmatrix}}
\newcommand{\evmx}{\end{vmatrix}}

\newcommand{\wh}{\widehat}
\newcommand{\wt}{\widetilde}
\newcommand{\f}{\frac}
\newcommand{\df}{\dfrac}

\newcommand{\imp}{\Longrightarrow}

\newcommand{\inc}{\subseteq}

\newcommand{\setm}{\setminus}

\newcommand{\Id}{\mathrm{Id}}

\newcommand{\minimize}[1]{\underset{#1}{\rm minimize}\,}

\newcommand{\dist}{\operatorname{dist}}

\newcommand{\La}{\Lambda}
\newcommand{\eps}{\varepsilon}
  
\newcommand{\ival}[2]{[{#1} \colon \hspace{-1mm} {#2} ]}

\title{\vspace{-20mm}Learning the Maximum of a H\"older Function from Inexact Data \vspace{-3mm}\\  \rule{6.5in}{1.2pt} \vspace{-15mm}}
\author{Simon Foucart\footnote{S. F. is partially supported by grants from the NSF (DMS-2505204) and the ONR (N00014-20-1-2787).}  --- Texas A\&M University}
\date{\vspace{-6mm}
\rule{100mm}{0.8pt}}

\newcommand\shorttitle{Learning the Maximum of a H\"older Function from Inexact Data}
\newcommand\authors{S. Foucart}

\newcommand{\rev}[1]{{#1}}

\begin{document}
\maketitle

%% Add abstract, keywords, and AMS classification
\vspace{-15mm}
\begin{abstract}
Within the theoretical framework of Optimal Recovery,
one determines in this article the {\em best} procedures to learn a quantity of interest depending on a H\"older function acquired via inexact point evaluations at fixed datasites.
{\em Best} here refers to procedures minimizing worst-case errors.
The elementary arguments hint at the possibility of tackling nonlinear quantities of interest,
with a particular focus on the function maximum.
In a local setting,
i.e., for a fixed data vector,
the optimal procedure
(outputting the so-called Chebyshev center)
is precisely described relatively to a general model of inexact evaluations.
Relatively to a slightly more restricted model and
in a global setting,
i.e., uniformly over all data vectors,
another optimal procedure is put forward,
showing how to correct the natural underestimate that simply returns the data vector maximum.
Jitterred data are also briefly discussed as a side product of evaluating the minimal worst-case error optimized over the datasites.
\end{abstract}



\vspace{-5mm}
\begin{center}
\rule{100mm}{0.8pt}
\end{center}

%%%%%%%%%%%%%%%%%
%% The main text starts here %%
%%%%%%%%%%%%%%%%%

\section{Introduction}

How can one reliably estimate the maximum $\max[f] := \max \{ f(x), x \in \cX \}$ of an unknown function from point evaluation data $y_m = f(x^{(m)})$, $m \in \ival{1}{M}$?
This question has been considered e.g. in \cite{CDPW,RMB}
with the common intent of avoiding the construction of a full approximant $\wh{f}$ to $f$ followed by the computation of the maximum of $\wh{f}$.
Both articles also assumed that the datasites $x^{(m)}$ could be freely chosen and that the data $y_m$ were exact.
In contrast,
in this note, the datasites are prescribed and the data are inexact, i.e., of the form
\be
\label{DefY}
y_m = f(x^{(m)}) + e_m
\qquad \mbox{ with }\quad |e_m| \le \eps_m
\ee
for some known $\eps_1,\ldots,\eps_M \ge 0$
(sometimes all taken to equal the same $\eps$, but definitely not always).
One shortly writes $\vy = {\bm \La}_{\vx} f  + \ve$,
with ${\bm \La}_{\vx}$ being the linear map assigning $[g(x^{(1)});\cdots;g(x^{(M)})] \in \bR^M$ to a function $g$
and with $\ve \in  \cE_{\bm \eps} := [-\eps_1,\eps_1] \times \cdots \times [-\eps_M, \eps_M]$.
Clearly, estimating $\max[f]$ from this information alone is nonsensical without an assumption preventing wild behaviors of $f$ away from the datasites.
Such an assumption is here expressed as $f \in \cK$,
where the so-called model set $\cK$ is the set of Lipschitz functions defined on a compact metric space $\cX$,
or more generally the set of H\"older functions with parameters $\alpha \in (0,1]$ given by\footnote{The set $\cH_\alpha$ can also be interpreted as a set of Lipschitz functions with respect to a modified metric defined by ${\rm dist}_\alpha(x,x') = \dist(x,x')^\alpha$, $x,x' \in \cX$.} 
\be
\label{DefH}
\cH_{\alpha} := 
\{ f:\cX \to \bR \mbox{ such that }
|f(x)-f(x')| \le \dist(x,x')^\alpha
\mbox{ for all } x,x' \in  \cX  \}.
\ee
With a focus on worst-case scenarios for assessing a procedure $\Delta: \vy \mapsto z$ designed to estimate $\max[f]$,
or another quantity of interest $\Gamma(f) \in Z$,
one defines worst-case errors 
\vspace{-5mm}
\begin{itemize}
\item locally at $\vy$
via
\be
\label{DefLWCE}
{\rm lwce}_{\vy,{\bm \eps}}(z) := \sup \big\{ \big \|\Gamma(f) - z \big\|_Z: f \in \cH_\alpha, \vy - {\bm \La}_{\vx} f \in \cE_{\bm \eps} \big\},\vspace{-2mm}
\ee
\item or globally via
\be
\label{DefGWCE}
{\rm gwce}_{\bm \eps}(\Delta) : = \sup \big\{ \big\|\Gamma(f) - \Delta( {\bm \La}_{\vx} f + \ve ) \big\| :  f \in \cH_\alpha, \ve \in \cE_{\bm \eps} \big\}.\vspace{-5mm}
\ee
\end{itemize}
Note that the `local' problem of minimizing ${\rm lwce}_{\vy, {\bm \eps}}(z)$ over $z$ is harder than solving the `global' problem of minimizing ${\rm gwce}_{\bm \eps}(\Delta)$ over $\Delta$,
in the sense that solutions $\wh{z}_{\vy, {\bm \eps}}$ to the former directly provide a solution $\Delta: \vy \mapsto \wh{z}_{\vy,{\bm \eps}}$ to the latter---in passing, one should take notice of the identity ${\rm gwce}_{\bm \eps}(\Delta) = \sup \{ {\rm lwce}_{\vy, {\bm \eps}}(\Delta(\vy)): \vy \in {\bm \La}_{\vx}(\cH_\alpha) + \cE_{\bm \eps} \}$.
The framework just described is the one underlying the theories of Optimal Recovery and subsequently of Information-Based Complexity 
(where the optimal selection of ${\bm \La}_{\vx}$ is also sought).
These theories often assume that the model set is convex and symmetric and that $\Gamma$ is a linear map,
implying in several situations that a globally optimal map $\Delta$ can be chosen linear
(see \cite[Chapter 4]{NovWoz}, \cite[Chapters 9 and 10]{BookDS}, and some of the author's recent works),
hinting once again at the preferability of the global problem over the local one.
The case of a linear quantity of interest $\Gamma$ treated here provides yet another situation where global optimality can be achieved by a linear $\Delta$,
see Theorem \ref{ThmLinearQ}.
Still, the locally optimal $\Delta$ is actually not very burdensome either,
see Theorem \ref{ThmLoc}.
In the case $\Gamma = \Id$ with exact data and a model class of Lipschitz functions,
this result was already obtained in \cite{Bel},
where the inspiration for the McShane-extension-type functions $\ell_{\vy, {\bm \eps}}$ and $u_{\vy, {\bm \eps}}$ of \eqref{DefLy}--\eqref{DefUy} comes from.
The novelty of this note
\rev{is the realization that,
even with inexact data, 
local and global Optimal Recovery problems may also be completely solved for some nonlinear quantities of interest,
in particular for the maximum of a function.
The situation of inexact data and of linear quantities of interest had been studied in \cite{Pla},
while the situation of exact data and of the maximum as quantity of interest had been treated in \cite[Section 1.3.6]{Nov} following the works \cite{Suk71,Was},
but amalgamating the two situations required a fresh perspective.
}

%beyond the extension of the result to inexact data and H\"older rather than Lipschitz model sets,
%is the realization that local and global Optimal Recovery problems may also be completely solved for some nonlinear quantities of interest---a situation rarely, if ever, considered in Optimal Recovery---namely for the maximum of a function.
% as put forward at the start.

Here is a preview of this note's contribution about the case $\Gamma(f) = \max[f]$.
First,
Section \ref{SecLoc} solves the local problem by showing,
as a specific instance of Theorem \ref{ThmLoc}, that
$$
\Delta^{\rm loc}:
\vy \in \bR^m \mapsto \f{\max[\ell_{\vy,{\bm \eps}}] + \max[u_{\vy,{\bm \eps}}]}{2} \in \bR
$$
is a locally optimal procedure for the estimation of $\max[f]$,
where 
\begin{align}
\label{DefLy}
\ell_{\vy,{\bm \eps}}(x) & := \max_{m \in \ival{1}{M}}  \big( y_m - \eps_m - \dist(x,x^{(m)})^\alpha  \big),\\
\label{DefUy}
u_{\vy,{\bm \eps}}(x) & := \min_{m \in \ival{1}{M}}  \big( y_m + \eps_m + \dist(x,x^{(m)})^\alpha  \big).
\end{align}

One may be content with using $\Delta^{\rm loc}$ as a global solution---it is!
Still, Section \ref{SecGlo} re-solves the global problem,
under the proviso that the $\eps_m$'s are all equal to a common $\eps$,
by exhibiting a `simpler' globally optimal procedure.
Namely, Theorem \ref{ThmGloMax} shows that
$$
\Delta^{\rm glo}:
\vy \in \bR^m \mapsto \max_{m \in \ival{1}{M}} \big( y_m \big) + \rev{\f{1}{2}} \max[U],
$$ 
is a globally optimal procedure for the estimation of $\max[f]$,
where the function $U$ is nothing else than $u_{{\bm 0}, {\bm 0}}$, i.e., it is given by\footnote{For $\alpha = 1$ and $\dist(x,x') = \|x-x'\|_{\ell_2^d}$,
the quantity $\max[U]$ is sometimes called the fill distance of $\{x^{(1)},\ldots,x^{(M)} \}$ for $\cX$, see \cite{Wen}.}
\be
\label{DefUU}
U(x)  := \min_{m \in \ival{1}{M}}  \big( \dist(x,x^{(m)})^\alpha \big).
\ee
This result makes perfect sense in hindsight,
at least for exact data ($\eps=0$).
Indeed, a natural estimation of $\max[f]$ is $\max_m(y_m)$,
but this always underestimates $\max[f]$,
so a correction term---depending on the location of the $x^{(m)}$'s---should be added.
Note that the map $\Delta^{\rm glo}$,
unlike $\Delta^{\rm loc}$,
is independent of $\eps$,
but its global worst-case error does depend on it.

Section \ref{SecExp} concludes by 
considering the situation where the metric $\dist$ stems from a norm on $\bR^d$ and the set $\cX$ is the unit ball for this norm.
The order of the smallest  possible global worst-case error over the choice of $x^{(1)},\ldots,x^{(M)}$ is uncovered there.
Perhaps unsurprisingly,
it features $\eps$ and $M^{-\alpha/d}$, see Theorem \ref{ThmGenX}.
In the particular case of the $\ell_\infty$-norm on $\bR^d$,
for which $\cX = [-1,1]^d$,
one can even give an exact expression for this smallest possible global worst-case error.
Perhaps unsurprisingly again,
the best choice for the $x^{(m)}$'s turn out to be the centers of the hypercubic cells partitioning $[-1,1]^d$ in a regular grid,
see Theorem \ref{ThmXCube} and its proof.

%Note one significant difference between the linear and nonlinear cases:
%in the linear case, optimally recovering the whole object $f$ and then applying $Q$ yields a optimal estimation of $Q(f)$,
%but this is not true e.g if $Q = \max$...\TODO{to check all of this...}

\section{Locally Optimal Solutions}
\label{SecLoc}

In this section,
considering a fixed observation vector $\vy \in \bR^m$ with entries given, as in \eqref{DefY}, by
$y_m = f(x^{(m)}) + e_m$, $|e_m| \le \eps_m$, 
one aims at producing an estimation of $\max[f]$---or more generally of a quantity of interest $\Gamma(f) \in Z$---which is worst-case optimal at this particular $\vy$.
The optimality, relatively to the model set $\cH_\alpha$ of H\"older functions introduced in \eqref{DefH},
thus refers to the local sense as described in \eqref{DefLWCE}\footnote{\rev{`Locally optimal' maps are also called `strongly optimal' or `central' in Information-Based Complexity, see \cite{TWW}.}}.
For full generality (not used beyond intuitive cases),
it is assumed that the normed space~$Z$ into which $\Gamma$ maps 
is a Banach lattice.
An overview and more details about Banach lattices can be found e.g. in \cite[Chapters 1 and 12]{HGBS},
\rev{see also \cite{Mey}},
but essentially these are normed spaces equipped with an order relation, compatible with vector addition and scalar multiplication, such that one can consider the `maximum' of two elements from $Z$, 
and hence the `absolute value' of an element $z \in Z$ as defined by $|z| = \max\{ -z, z\} \in Z$.
The norm should be monotone\footnote{There are several notions of monotone norms, listed and compared in \cite{JohNyl}.}, in the sense that,
for any $z,z' \in Z$, 
\be
\label{MonoProp}
|z| \ge |z'| \imp \|z\|_Z \ge \|z'\|_Z.
\ee
As for the quantity of interest $\Gamma: B(\cX) \to Z$
with domain equal to the space of bounded functions on $\cX$,
it is not assumed to be linear.
Instead, it is only required to be monotone,
in the sense that, for any $f,f' \in B(\cX)$,
\be
\label{QMono}
f \ge f' \imp \Gamma(f) \ge \Gamma(f').
\ee
Note that this covers quite a number of cases:
as linear examples, 
one can think of $\Gamma$ being the identity (i.e., $\Gamma(f)=f$),
a domain restriction (i.e., $\Gamma(f)=f_{| \Omega}$),
and a positive linear functional (e.g. $\Gamma(f) = \int_\cX f$), etc.;
as nonlinear examples,
one can think of $\Gamma$ being an increasing map applied pointwise (e.g. $\Gamma(f) = \exp(f)$)
or---the main focus of this note---the maximum of a function (i.e., $\Gamma(f) = \max[f]$).
The result below exactly solves the problem of minimizing the local worst-case error.
The function $\ell_{\vy,{\bm \eps}}$ and $u_{\vy,{\bm \eps}}$ appearing there are the ones introduced in \eqref{DefLy}--\eqref{DefUy}.

\bthm
\label{ThmLoc}
For the recovery of a monotone quantity of interest $\Gamma$ mapping into a Banach lattice~$Z$,
given the data \eqref{DefY} and the model set \eqref{DefH},
a minimizer of ${\rm lwce}_{\vy,{\bm \eps}}(z)$ over all $z \in Z$ is given by
$$
\wh{z}_{\vy, {\bm \eps}} =  \df{\Gamma(\ell_{\vy,{\bm \eps}}) + \Gamma(u_{\vy,{\bm \eps}})}{2}, 
$$
while the value of the minimum is
$$ 
\wh{r}_{\vy, {\bm \eps}} = \f{\| \Gamma(u_{\vy,{\bm \eps}}) - \Gamma(\ell_{\vy,{\bm \eps}}) \|_Z}{2}.
$$
\ethm

\bpf
Central to the argument is a basic observation 
(elucidating the notation $\ell_{\vy,{\bm \eps}}$ and $u_{\vy,{\bm \eps}}$, which stands for lower and upper bounds, respectively).
It starts by writing, for any $f \in \cH_\alpha$, any $x \in \cX$,
and any $m \in \ival{1}{M}$,
$$
|f(x) - y_m | \le |f(x^{(m)}) - y_m| + |f(x) - f(x^{(m)})|
\le \eps_m + \dist(x,x^{(m)})^\alpha ,
$$
that is to say
$$
 y_m - \eps_m - \dist(x,x^{(m)})^\alpha 
\le f(x) 
\le y_m + \eps_m + \dist(x,x^{(m)})^\alpha .
$$
Taking the maximum and minimum over $m \in \ival{1}{M}$,
one obtains
\be
\label{lb<f<ub}
\ell_{\vy,{\bm \eps}}(x) \le f(x) \le u_{\vy,{\bm \eps}}(x)
\qquad \mbox{for all } f \in \cH_\alpha \mbox{ and } x \in \cX.
\ee
Furthermore, one observes that 
both $\ell_{\vy,{\bm \eps}}$ and $u_{\vy,{\bm \eps}}$ are data-consistent and model-consistent,
i.e.,
$|y_m - \ell_{\vy, {\bm \eps}}(x^{(m)})| \le \eps_m$ and $|y_m - u_{\vy, {\bm \eps}}(x^{(m)})| \le \eps_m$ for all $m \in \ival{1}{M}$,
as well as 
 $\ell_{\vy,{\bm \eps}} \in \cH_\alpha$ and $u_{\vy,{\bm \eps}} \in \cH_\alpha$. 
Indeed, for the first statement,
using the defining expressions of $\ell_{\vy,{\bm \eps}}$ and $u_{\vy,{\bm \eps}}$ evaluated at $x^{(m)}$ for any $m \in \ival{1}{M}$,
one sees that $y_{m} - \eps_m \le \ell_{\vy,{\bm \eps}}(x^{(m)}) $ 
and $u_{\vy,{\bm \eps}}(x^{(m)}) \le y_m + \eps_m$,
so that, in view of \eqref{lb<f<ub},
one arrives at $y_m - \eps_m \le \ell_{\vy,{\bm \eps}}(x^{(m)}) \le u_{\vy,{\bm \eps}}(x^{(m)}) \le y_m + \eps_m$ for any $m \in \ival{1}{M}$. This is data-consistency.
For the second statement,
given $x,x' \in \cX$ and  $m \in \ival{1}{M}$,
one has
\begin{align*}
\big( y_m - \eps_m - \dist(x,x^{(m)})^\alpha  \big) - \ell_{\vy,{\bm \eps}}(x')
& \le \big( y_m - \eps_m - \dist(x,x^{(m)})^\alpha  \big) - \big( y_m - \eps_m - \dist(x',x^{(m)})^\alpha  \big)\\
& = -\dist(x,x^{(m)})^\alpha + \dist(x',x^{(m)})^\alpha \le  \dist(x,x')^\alpha,
\end{align*}
which yields
$\ell_{\vy,{\bm \eps}}(x) - \ell_{\vy,{\bm \eps}}(x') \le \dist(x,x')^\alpha$
 by taking the maximum over $m \in \ival{1}{M}$.
 The inequality $| \ell_{\vy,{\bm \eps}}(x) - \ell_{\vy,{\bm \eps}}(x') |  \le \dist(x,x')^\alpha$ is obtained by exchanging the roles of $x$ and $x'$,
thus showing that $\ell_{\vy,{\bm \eps}} \in \cH_{\alpha}$.
This is model-consistency for $\ell_{\vy,{\bm \eps}}$.
The case of $u_{\vy,{\bm \eps}}$ is treated similarly.

Armed with the above observations, 
it is time to prove the prospective result that, for any $z \in Z$,
$$
{\rm lwce}_{\vy,{\bm \eps}}(z) \ge \f{\|\Gamma(u_{\vy,{\bm \eps}}) - \Gamma(\ell_{\vy,{\bm \eps}})\|_Z}{2}
\qquad \mbox{with equality when }
z = \wh{z}_{\vy,{\bm \eps}} = \f{\Gamma(\ell_{\vy,{\bm \eps}}) + \Gamma(u_{\vy,{\bm \eps}})}{2}. 
$$ 
To justify one part of this claim, given any $z \in Z$,
the fact that both $\ell_{\vy,{\bm \eps}}$ and $u_{\vy,{\bm \eps}}$ are model- and data-consistent implies that
$$
{\rm lwce}_{\vy,{\bm \eps}}(z)
\ge \begin{cases}
\|\Gamma(u_{\vy,{\bm \eps}}) - z\|_Z \\ \|\,\Gamma(\ell_{\vy,{\bm \eps}}) - z\|_Z
\end{cases} \hspace{-3mm}
\ge \f{1}{2} \big( \|\Gamma(u_{\vy,{\bm \eps}}) - z\|_Z + \|\Gamma(\ell_{\vy,{\bm \eps}}) - z\|_Z \big)
\ge \f{1}{2} \|\Gamma(u_{\vy,{\bm \eps}}) - \Gamma(\ell_{\vy,{\bm \eps}})\|_Z.
$$
For the other part of the claim,
given any $f \in \cH_\alpha$ satisfying $\vy - {\bm \La}_{\vx} f \in \cE_{\bm \eps}$,
i.e.,
$|y_m - f(x^{(m)})| \le \eps_m$ for all $m \in \ival{1}{M}$,
the inequalities \eqref{lb<f<ub}
and the monotonicity \eqref{QMono} of $\Gamma$
yield
$\Gamma(\ell_{\vy,{\bm \eps}}) \le \Gamma(f) \le \Gamma(u_{\vy,{\bm \eps}})$,
i.e.,
$$
- \f{\Gamma(u_{\vy,{\bm \eps}})-\Gamma(\ell_{\vy,{\bm \eps}})}{2}
\le \Gamma(f) - \f{\Gamma(\ell_{\vy,{\bm \eps}}) + \Gamma(u_{\vy,{\bm \eps}})}{2}
\le \f{\Gamma(u_{\vy,{\bm \eps}})-\Gamma(\ell_{\vy,{\bm \eps}})}{2}.
$$
From here, it follows (with the `absolute value' being the one relative to the Banach lattice $Z$) that
$$
\left| \Gamma(f) - \f{\Gamma(\ell_{\vy,{\bm \eps}}) + \Gamma(u_{\vy,{\bm \eps}})}{2} \right|
\le \f{\Gamma(u_{\vy,{\bm \eps}})-\Gamma(\ell_{\vy,{\bm \eps}})}{2},
$$
allowing to conclude, using the monotonicity \eqref{MonoProp} of the norm on $Z$, that 
$$
\left\| \Gamma(f) - \f{\Gamma(\ell_{\vy,{\bm \eps}}) + \Gamma(u_{\vy,{\bm \eps}})}{2} \right\|_Z
\le \f{\| \Gamma(u_{\vy,{\bm \eps}})-\Gamma(\ell_{\vy,{\bm \eps}}) \|_Z}{2}.
$$
Taking the supremum over $f$ shows that 
${\rm lwce}_{\vy,{\bm \eps}}(\wh{z}_{\vy, {\bm \eps}}) \le \|\Gamma(u_{\vy,{\bm \eps}}) - \Gamma(\ell_{\vy,{\bm \eps}})\|_Z/2$,
thus justifying the other part of the claim.
The proof is now complete.
\epf

The elements $\wh{z}_{\vy, {\bm \eps}} \in Z$ and $\wh{r}_{\vy,{\bm \eps}} := {\rm lwce}_{\vy,{\bm \eps}}(\wh{z}_{\vy, {\bm \eps}})  \in \bR_+$
are called Chebyshev center and radius, respectively, for the set $\{ \Gamma(f): f \in \cH_\alpha, \, \vy - {\bm \La}_{\vx} f \in \cE_{\bm \eps} \}$,
as they are solutions to
$$
\minimize{z,r} \; r
\quad
\mbox{subject to } \|\Gamma(f) - z\|_Z \le r \mbox{ whenever }
f \in \cH_\alpha \mbox{ and } \vy - {\bm \La}_{\vx} f \in \cE_{\bm \eps} .
$$
In other words, they are center and radius of a smallest $Z$-ball containing this set.
Using the notation $\wh{f}_{\vy, {\bm \eps}}$
for the Chebyshev center in the special case $\Gamma = I$,
Theorem \ref{ThmLoc} reveals in particular that,
if the quantity of interest $\Gamma$ is linear,
then $\wh{z}_{\vy, {\bm \eps}} = \Gamma(\wh{f}_{\vy, {\bm \eps}})$,
i.e., a locally optimal estimation of $\Gamma(f)$ can be obtained by construction a full approximant and plugging it into $\Gamma$.
This fortuity does not hold anymore for nonlinear quantities of interest, say for $\Gamma = \max$.
Indeed, the plug-in estimate $\max[\wh{f}_{\vy, {\bm \eps}}] = \max[(\ell_{\vy, {\bm \eps}} + u_{\vy, {\bm \eps}})/2]$
is generally smaller than $\wh{z}_{\vy, {\bm \eps}} = (\max[\ell_{\vy, {\bm \eps}}] + \max[u_{\vy, {\bm \eps}}] )/2$,
as confirmed from Figure \ref{fig} below. \vspace{-5mm}
   
\begin{figure}[!htbp]
\centering     
\subfigure[Full recovery of $f$]{\includegraphics[width=80mm]{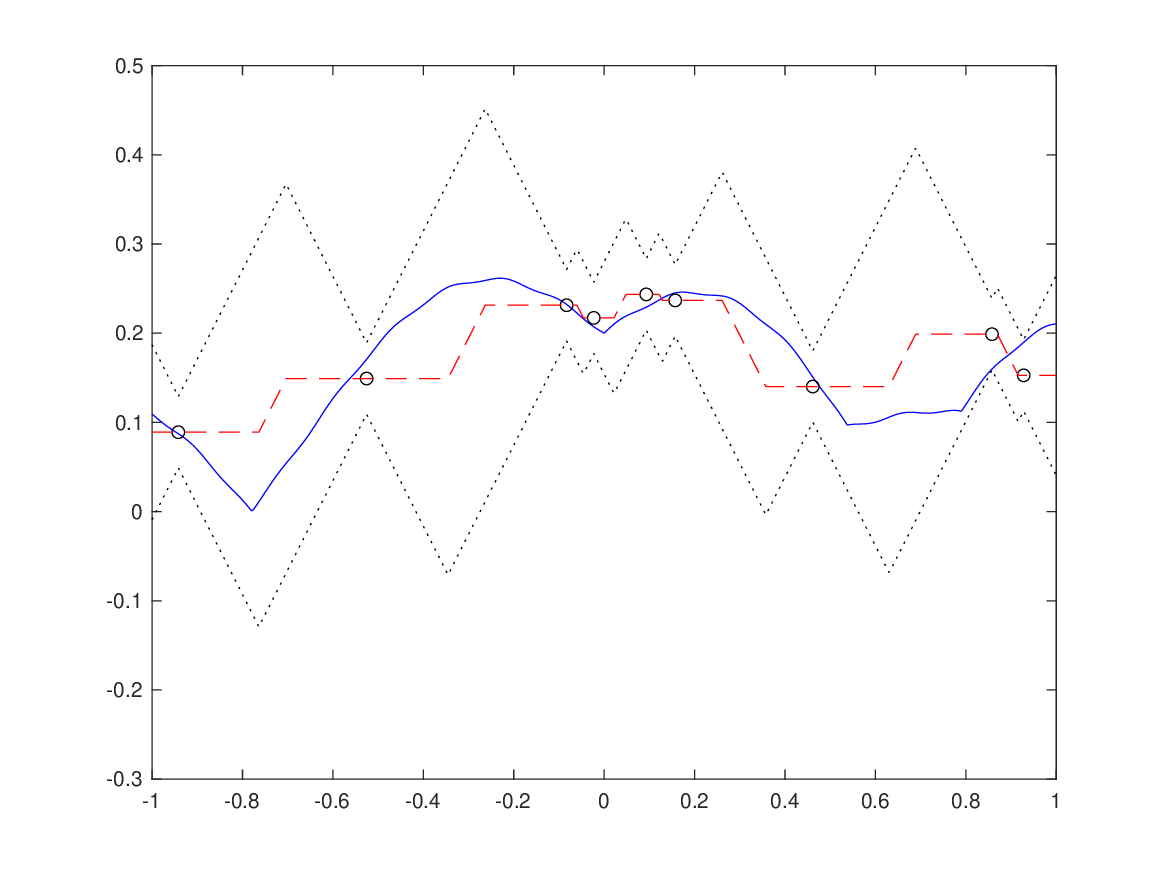}}
\subfigure[Estimation of $\max{[}f{]} $]{\includegraphics[width=80mm]{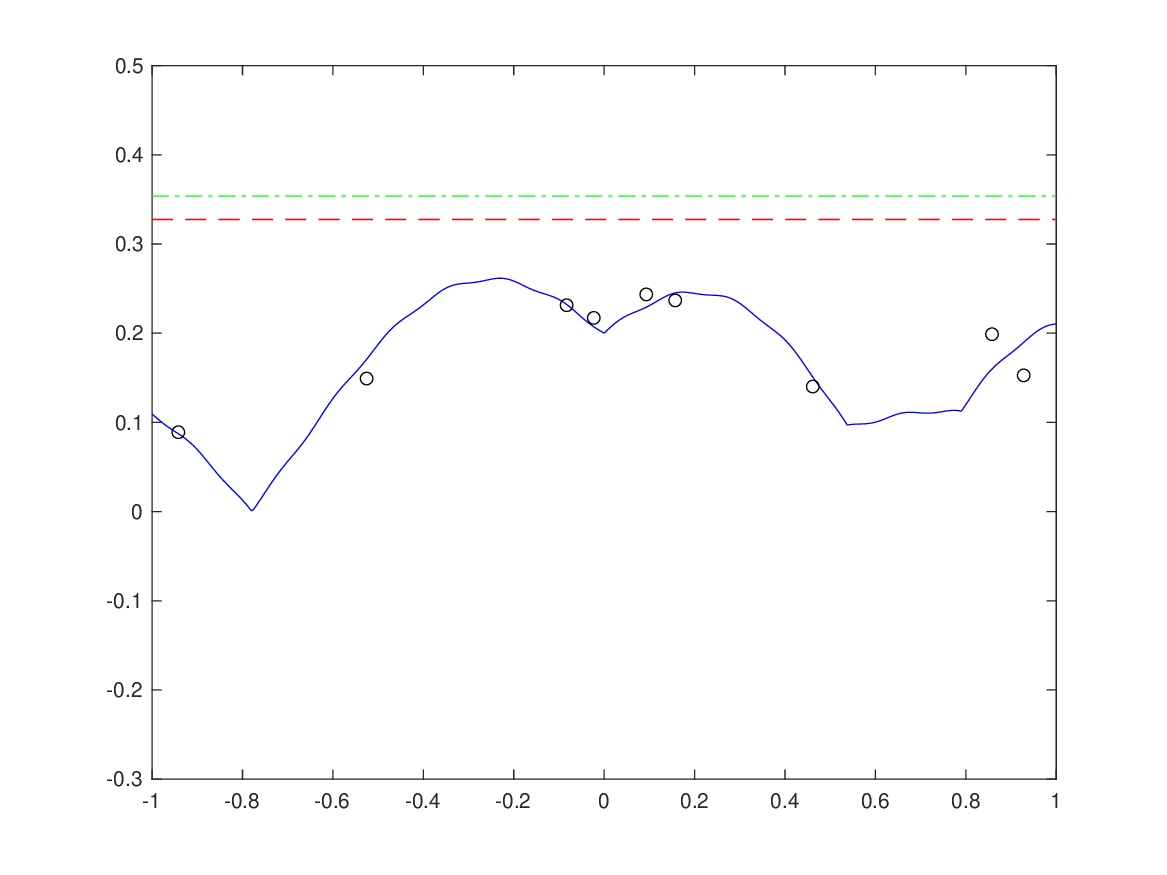}}
\caption{An example of a univariate Lipschitz function $f$ (solid blue curve) acquired via inexact point values (black circles).
The left panel displays the lower function $\ell_{\vy,{\bm \eps}}$ and the upper function $u_{\vy,{\bm \eps}}$ from \eqref{DefLy}-\eqref{DefUy} (dotted curves),
as well as the locally optimal estimation of $f$ from Theorem~\ref{ThmLoc} (dashed red curve).
The right panel shows, as horizontal lines, 
the locally optimal estimation of $\max[f]$ from Theorem~\ref{ThmLoc} (dashed red line)
and the globally optimal estimation of $\max[f]$ from Theorem~\ref{ThmGloMax} (dash-dotted green line).}
\label{fig}
\end{figure}

%\brk
%For later purposes,
%one notices---again in the noiseless case---that
%$$
%\max \Big\{ \f{\ell_y + u_y}{2} \Big\} = \max_m\{y_m\}.
%$$ 
%To see this,
%for $x \in \cX$,
%consider an index $m$ achieving the minimum in the defining expression of $u_y(x)$,
%one has
%$$
%\f{1}{2} \big( \ell_y(x)+u_y(x) \big)
%\le \f{1}{2} \big( (y_m - \eps_m - \dist(x,x^{(m)})) + (y_m + \eps_m + \dist(x,x^{(m)}) \big) 
%= y_m \le \max_m\{y_m\}.
%$$ 
%The inequality $\max{(\ell_y + u_y)/2} \le \max_m \{y_m\}$ follows.\TODO{this was valid in the noisy setting, too.}
%
%For the reverse inequality,
%consider an index $m$ such that $y_m = \max_m\{ y_m \}$.
%Recall that the inequality $\ell_y(x) \le u_y(x)$ implies that
%$y_m - \eps_m - \dist(x,x^{(m)}) \le y_{m'} + \eps_{m'} + \dist(x,x^{(m')})$,
%which implies---in the noiseless case---
%by taking $x=x^{(m)}$, that 
%$y_m \le y_{m'} + \dist(x^{(m)},x^{(m')})$ for all $m'$,
%and hence $u_y(x^{(m)}) = y_m$.
%As a result,
%$$
%\max \Big\{ \f{\ell_y + u_y}{2} \Big\}
%\ge  \f{\ell_y(x^{(m)}) + u_y(x^{(m)})}{2}
%\ge \f{y_m+y_m}{2} = y_m = \max_m \{ y_m\}.
%$$
%\TODO{Easier: from \eqref{lb<f<ub} and the defining expressions of $\ell_y$ and $u_y$, we have 
%$y_m - \eps_m \le \ell_y(x^{(m)}) \le u_y(x^{(m)}) \le y_m + \eps_m$...}
%\erk

The Chebyshev center $\wh{z}_{\vy, {\bm \eps}}$ is not supplied as explicitly as ideally desired,
even for $\Gamma(f) = \max[f]$.
For the specific case of univariate ($\cX = [-1,1]$)
Lipschitz  ($\alpha=1$) functions  acquired via exact point evaluations ($\eps_m = 0$), 
it is actually possible to be completely explicit.

\bprop
\label{PropUnivariate}
For functions $f:[-1,1] \to \bR$
satisfying the condition $|f(x)-f(x')| \le |x-x'|$ for all $x,x' \in [-1,1]$
and observed via $y_1=f(x^{(1)})$, ..., $y_M=f(x^{(M)})$
at $x^{(1)} < \cdots < x^{(M)} \in (-1,1)$,
the locally optimal estimation of $\max[f]$ is expressed as
$$
\wh{z}_{\vy, {\bm 0}} = \f{1}{2}
\Big[
\max_{m \in \ival{1}{M}} \big( y_m \big)
+ \max_{m \in \ival{0}{M}} \Big( \f{x^{(m+1)} - x^{(m)} + y_{m+1} + y_m}{2} \Big)
\Big],
$$
having set $x^{(0)}=-1$, $y_0 = 0$, $x^{(M+1)}=1$, and $y_{M+1}=0$. 
\eprop 

\bpf
Based on Theorem \ref{ThmLoc},
one needs to untangle $\max[\ell_{\vy, {\bm 0}}]$
and $\max[u_{\vy, {\bm 0}}]$.
The former can be dealt with in full generality.
Indeed, in view of $\ell_{\vy, {\bm 0}}(x) = \max_m \big( y_m - \dist(x,x^{(m)})^\alpha \big)$,
one observes on the one hand that $\ell_{\vy, {\bm 0}}(x) \le \max_m \big( y_m  \big)$ for all $x \in \cX$ 
and on the other hand that $\ell_{\vy, {\bm 0}}(x^{(m)}) \ge y_m $ for all $m \in \ival{1}{M}$.
These two facts combine to show that 
$$
\max[\ell_{\vy, {\bm 0}}] = \max_m \big( y_m \big) .
$$
It remains to determine $\max[u_{\vy, {\bm 0}}]$ explicitly
in the specific univariate setting of the proposition.
Precisely,
one shall show that 
\be
%\label{Need1}
\label{Need}
\max[u_{\vy, {\bm 0}}]
%= \max_{x \in [-1,1]} \min_{m \in \ival{1}{M}} \big( y_m + |x-x^{(m)}| \big)
%\label{Need2}
= \max_{m \in \ival{0}{M}} \Big( \f{x^{(m+1)} - x^{(m)} + y_{m+1} + y_m}{2} \Big).
\ee
To this end, one considers  
the Voronoi-type cells $V_{\vy,1},\ldots,V_{\vy, M}$ defined
for $m \in \ival{1}{M}$ by
$$
V_{\vy,m} := \big\{ x \in [-1,1]:
y_m + |x - x^{(m)} | 
\le y_n + |x-x^{(n)}|
\; \mbox{ for all } n \in \ival{1}{M} 
\big\},
$$
so that 
\be
\label{MaxViaVor}
\max[u_{\vy, {\bm 0}}]
= \max_{x \in [-1,1]} \min_{m \in \ival{1}{M}} \big( y_m + |x-x^{(m)}| \big)
= \max_{m \in \ival{1}{M}} \max_{x \in V_{\vy,m}} \big( y_m + |x-x^{(m)}| \big).
\ee
At a fixed $m \in \ival{1}{M}$,
introducing $V_{\vy,m,n}:= \{ x \in [-1,1]: |x-x^{(m)}| - |x - x^{(n)}| \le y_n - y_m \}$
for each $n \in \ival{1}{M}$,
one has $V_{\vy,m} = \cap_{n \not= m} V_{\vy,m,n}$.
For $n<m$,
since $x \mapsto |x-x^{(m)}| - |x - x^{(n)}|$ is a ramp function joining $|x^{(n)}-x^{(m)}|$ on the left of $x^{(n)}$
to $-|x^{(n)}-x^{(m)}|$ on the right of $x^{(m)}$
and since $|y_n - y_m| \le |x^{(n)} - x^{(m)}|$ by the Lipschitz condition,
a point $x \in [-1,1]$ belongs to $V_{\vy, m,n}$
if and only if $x \ge \alpha^{(m,n)}$,
where  $\alpha^{(m,n)} \in [x^{(n)},x^{(m)}]$ is implicitly defined by the equation $|\alpha^{(m,n)}-x^{(m)}| - |\alpha^{(m,n)} - x^{(n)}| = y_n - y_m$,
or explicitly given by $\alpha^{(m,n)} = (x^{(m)}+x^{(n)} + y_m-y_n)/2$.
For $n>m$, one also obtains that a point $x \in [-1,1]$ belongs to $V_{\vy, m,n}$
if and only if $x \le \beta^{(m,n)}$,
where 
$\beta^{(m,n)} \in [x^{(m)},x^{(n)}]$ is implicitly defined by the equation 
$|\beta^{(m,n)}-x^{(m)}| - |\beta^{(m,n)} - x^{(n)}| = y_n - y_m$,
or explicitly given by 
$\beta^{(m,n)} = (x^{(m)}+x^{(n)} - y_m + y_n)/2$.
As a result, for $m \in [2:M-1]$, one derives that
\begin{align*}
V_{\vy,m} & = \Big( \bigcap_{n < m} [\alpha^{(m,n)},1] \Big)\bigcap \Big( \bigcap_{n>m} [-1,\beta^{(m,n)}] \Big)
= \big[ \max_{n<m} \alpha^{(m,n)},1 \big] \bigcap \big[-1, \min_{n>m} \beta^{(m,n)} \big] \\
& 
% = [\alpha^{(m,m-1)},1] \bigcap [-1,\beta^{(m,m+1)}],
= [\alpha^{(m,m-1)},\beta^{(m,m+1)}],
\end{align*}
where one used the easily verifiable facts that $\alpha^{(m,n)} \le \alpha^{(m,m-1)}$ for $n<m$
and $\beta^{(m,n)} \ge \beta^{(m,m+1)}$ for $n>m$.
Now, noticing that $\beta^{(m,m+1)}$ coincides with $\alpha^{(m+1,m)}$ with equal value
$$
\xi^{(m)} := \f{x^{(m)}+x^{(m+1)}-y_m+y_{m+1}}{2},
\qquad m \in \ival{1}{M-1}, 
$$
one deduces that the Voronoi-type cells reduce to 
$$
V_{\vy, m} = [\xi^{(m-1)},\xi^{(m)}]
\qquad \mbox{for } m \in [2:M-1].
$$
This remains valid for $m=1$ if one sets $\xi^{(0)}=-1$
and for $m=M$ if one sets $\xi^{(M)}=1$---for instance,  $V_{\vy,M} = \cap_{n<M} [\alpha^{(M,n)},1] = [\alpha^{(M,M-1)},1] = [\xi^{(M-1)},1]$.
It follows that,
for any $m \in \ival{1}{M}$,
\begin{align*}
\max_{x \in V_{\vy,m}} \big( y_m + |x-x^{(m)}| \big)
& = \max\{ y_m + x^{(m)} - \xi^{(m-1)}, y_m + \xi^{(m)} - x^{(m)} \}\\
& = \max \Big\{ 
\f{x^{(m)} - x^{(m-1)} + y_{m} + y_{m-1}}{2},
\f{x^{(m+1)} - x^{(m)} + y_{m+1} + y_m}{2}
 \Big\},
\end{align*}
where the boundary cases $m=1$ and $m=M$ require some special attention.  
Substituting the latter into \eqref{MaxViaVor} and 
eliminating the repetitions leads to the expression announced in \eqref{Need}.
\epf

\section{Globally Optimal Solutions}
\label{SecGlo}

In this section,
one still aims at producing a recovery map $\Delta: \bR^m \to Z$
which is worst-case optimal for the estimation of a quantity of interest $\Gamma(f) \in Z$,
but the optimality  here refers to the global sense as described in \eqref{DefGWCE}.
As mentioned earlier,
the locally optimal recovery map $\Delta^{\rm loc}: \vy \in \bR^m \mapsto \wh{z}_{\vy, {\bm \eps}} \in Z$ from Theorem \ref{ThmLoc} is already automatically globally optimal.
Nevertheless, one may be interested in other maps $\Delta: \bR^m \to Z$ minimizing ${\rm gwce}_{\bm \eps}(\Delta)$,
one reason being that they could be easier to construct and manipulate,
e.g. they could be linear in selected instances where the quantity of interest $\Gamma$ is itself linear.
Contrary to Section~\ref{SecLoc},
the cases of linear and nonlinear $\Gamma$
can seemingly not be treated in one swift stroke,
so they are separated from now on.

\subsection{The case of linear quantities of interest}

In the result below,
the generic function $\mathbbm{1}_{V}$ denotes the indicator function of a set $V$,
meaning that $\mathbbm{1}_{V}(x) = 1$ if $x \in V$ and $\mathbbm{1}_{V}(x) = 0$ if $x \not\in V$,
while the sets $V_{{\bm \eps},1},\ldots,V_{{\bm \eps},M} \inc \cX$ are the Voronoi-type cells defined by
$$
V_{{\bm \eps},m} := \big\{ x \in \cX:
\eps_m + \dist(x,x^{(m)})^\alpha  
\le \eps_n + \dist(x,x^{(n)})^\alpha
\; \mbox{ for all } n \in \ival{1}{M} 
\big\}.
$$

\bthm
\label{ThmLinearQ}
For the recovery of a monotone linear quantity of interest $\Gamma$ from $B(\cX)$ into a Banach lattice $Z$,
given the datasites $x^{(1)},\ldots,x^{(M)} \in \cX$ and the model set \eqref{DefH},
the  minimal global worst-case error is 
$$
\inf_{\Delta: \bR^M \to Z}
{\rm gwce}_{{\bm \eps}}(\Delta)  = \|\Gamma(u_{{\bm 0},{\bm \eps}})\|_Z,
\qquad \mbox{ where }
u_{{\bm 0},{\bm \eps}}(x) = \min_{m \in \ival{1}{M}} \big( \eps_m + \dist(x, x^{(m)})^\alpha  \big).
$$
Furthermore, the infimum is achieved for the recovery map
$\Gamma \circ \Delta^{\rm glo}: \bR^M \to Z$, where
$$
\Delta^{\rm glo} : \vy \in \bR^M \mapsto \sum_{m=1}^M y_m \mathbbm{1}_{V_{{\bm \eps},m}} \in B(\cX).
$$
\ethm

\bpf
For the first part of the proof,
i.e., the fact that ${\rm gwce}_{{\bm \eps}}(\Delta) \ge \|\Gamma(u_{{\bm 0},{\bm \eps}})\|_Z$ for any $\Delta: \bR^m \to Z$,
one relies of the aforementioned identity 
${\rm gwce}_{\bm \eps}(\Delta) = \sup \{ {\rm lwce}_{\vy, {\bm \eps}}(\Delta(\vy)): \vy \in {\bm \La}_{\vx}(\cH_\alpha) + \cE_{\bm \eps} \}$.
One immediately derives that
$$
{\rm gwce}_{\bm \eps}(\Delta)
\ge {\rm lwce}_{{\bm 0}, {\bm \eps}}(\Delta(\vy))
\ge \wh{r}_{{\bm 0}, {\bm \eps}} = 
\f{\| \Gamma(u_{{\bm 0},{\bm \eps}}) - \Gamma(\ell_{{\bm 0},{\bm \eps}}) \|_Z}{2}
= \| \Gamma(u_{{\bm 0},{\bm \eps}})  \|_Z,
$$
where the last equality used  the fact (seen from \eqref{DefLy}--\eqref{DefUy}) that $\ell_{{\bm 0},{\bm \eps}} = - u_{{\bm 0},{\bm \eps}}$ and the linearity of $\Gamma$.

For the second part of the proof,
it is enough to show that ${\rm gwce}_{\bm \eps}(\Gamma \circ \Delta^{\rm glo}) \le \|\Gamma(u_{{\bm 0},{\bm \eps}})\|_Z$.
To this end,
let $f \in \cH_\alpha$, $\ve \in \cE_{\bm \eps}$,
and let the shorthand $\vy := {\bm \La}_{\vx} f + \ve$ be employed below.
Then, for any $x \in \cX$,
select $n \in \ival{1}{M}$ such that $x \in V_{{\bm \eps},n}$
and observe that
\begin{align*}
| f(x) - \Delta^{\rm glo}(\vy)(x) | & = 
|f(x) - \sum_{m=1}^M y_{m} \mathbbm{1}_{V_{{\bm \eps},m}}(x)|
= |f(x) - y_n| 
= |f(x)- (f(x^{(n)})+ e_n)|\\
& \le |e_n| + |f(x) - f(x^{(n)})| 
\le \eps_n + \dist(x,x^{(n)})^\alpha 
= u_{{\bm 0},{\bm \eps}}(x).
\end{align*}
This reads $-u_{{\bm 0},{\bm \eps}} \le f - \Delta^{\rm glo}(\vy) \le u_{{\bm 0},{\bm \eps}}$
and, thanks the monotonicity \eqref{QMono} and linearity of $\Gamma$,
it follows that $-\Gamma(u_{{\bm 0},{\bm \eps}}) \le \Gamma(f) - (\Gamma \circ \Delta^{\rm glo})(\vy) \le \Gamma(u_{{\bm 0},{\bm \eps}})$,
i.e., $| \Gamma(f) - (\Gamma \circ \Delta^{\rm glo})(\vy) | \le \Gamma(u_{{\bm 0},{\bm \eps}})$.
Finally,
thanks to the monotonicity \eqref{MonoProp} of the norm,
one concludes that $\| \Gamma(f) - (\Gamma \circ \Delta^{\rm glo})(\vy) \|_Z \le \| \Gamma(u_{{\bm 0},{\bm \eps}}) \|_Z$.
Taking the maximum over $f$ and $\ve$ finishes the proof.
\epf

A few comments are appropriate to conclude on the subject of monotone linear quantities of interest:\vspace{-5mm}
\begin{itemize}
\item the two optimal recovery maps that were exposed,
the locally---hence globally---optimal one from Theorem \ref{ThmLoc}
and the globally optimal one from Theorem \ref{ThmLinearQ},
both take the form $\Gamma \circ \Delta$, with $\Delta$ given by either
$$
\Delta^{\rm loc}(\vy) = \f{\ell_{\vy,{\bm \eps}} + u_{\vy,{\bm \eps}}}{2}
\qquad \mbox{ or } \qquad
\Delta^{\rm glo}(\vy) = \sum_{m=1}^M y_m \mathbbm{1}_{V_{{\bm \eps},m}},
$$
i.e., they are obtained by applying the quantity of interest after producing a full estimation;\vspace{-2mm}
\item when all the $\eps_m$'s are equal to a common $\eps$,
both $\Delta^{\rm loc}$ and $\Delta^{\rm glo}$ are independent of this $\eps$ (but of course the minimal worst-case errors do depend on $\eps$);\vspace{-2mm}
\item the local estimate $\Delta^{\rm loc}(\vy)$ is data- and model-consistent,
while the global estimate $\Delta^{\rm glo}(\vy)$ is not model-consistent, as it is discontinuous, hence not H\"older;\vspace{-2mm}
\item the global estimate $\Delta^{\rm glo}(\vy)$ depends linearly on $\vy \in \bR^m$,
while local estimate $\Delta^{\rm loc}(\vy)$ does not---this does not necessarily mean, however, that the former is simpler than the latter,
as producing $\Delta^{\rm loc}(\vy)$ via the defining expressions \eqref{DefLy}--\eqref{DefUy} of $\ell_{\vy,{\bm \eps}}$ and $u_{\vy,{\bm \eps}}$ can be viewed as easier than producing $\Delta^{\rm glo}(\vy)$, if done through the construction all the Voronoi-type cells;\vspace{-2mm}
\item the local estimate $\Delta^{\rm loc}(\vy)$ seems more convenient for streaming data
than the global estimate $\Delta^{\rm glo}(\vy)$,
as the former is easy to update (via the defining expressions \eqref{DefLy}--\eqref{DefUy} of $\ell_{\vy,{\bm \eps}}$ and $u_{\vy,{\bm \eps}}$)
when a new $(x^{(M+1)},y_{M+1})$ is added,
while this addition  potentially affects a large number of Voronoi-type cells for the latter.
But if one is primarily interested in prediction at a fixed $x \in \cX$,
then the global estimate is also straightforward to update,
as it essentially returns the value $y_m$ at the datasite $x^{(m)}$ closest to $x$,
so one only needs to check if $x$ is closer to $x^{(M+1)}$ than to the previous closest datasite.
\end{itemize}

\subsection{The case of the maximum}

Turning to the estimation of nonlinear quantities of interest,
as an encompassing result seems out of reach,
one concentrates  specifically on the maximum,
i.e., $\Gamma(f)=\max[f]$.
To start, consider the situation of exact observations,
i.e., $y_m = f(x^{(m)})$ for all $m \in \ival{1}{M}$.
To estimate $\max[f]$, 
it might seem natural to simply output $\max_{m} ( y_m )$,
but this is always an underestimation.
One probably needs to add a slight correction.
To guess what this correction could be,
look at  the local optimality result instantiated at $\vy={\bm 0}$,
which suggests $(\max[\ell_{{\bm 0},{\bm 0}}] + \max[u_{{\bm 0},{\bm 0}}])/2$.
Since $\max [ \ell_{{\bm 0},{\bm 0}} ] = 0$ and $u_{{\bm 0},{\bm 0}}$ coincides with the function $U$ defined in \eqref{DefUU},
this correction simplifies to $\max[U]/2$.
It turns to be the right guess,
remaining valid when the observation error bounds $\eps_m$ are all equal to a common~$\eps$.
In this case, 
a globally optimal recovery map can be chosen independently of this $\eps$, as formally stated below.
% (but of course the minimal global worst-case error does depend on $\eps$).

\bthm
\label{ThmGloMax}
For the recovery of the quantity of interest $f \in B(\cX) \mapsto \max[f] \in \bR$,
given the data~\eqref{DefY} with $\eps_m = \eps$ for all $m \in \ival{1}{M}$ and the model set \eqref{DefH},
the minimal global worst-case error is 
$$
\inf_{\Delta: \bR^M \to \bR}{\rm gwce}_{{\bm \eps}}(\Delta)  
= \f{1}{2} \max[U] + \eps,
\qquad \mbox{ where }
U(x) = \min_{m \in \ival{1}{M}} \big( \dist(x, x^{(m)})^\alpha \big).
$$
Furthermore, the infimum is achieved for the recovery map
$$
\Delta^{\rm glo} : \vy \in \bR^M \mapsto
\max_{m \in \ival{1}{M}} \big( y_m \big) + \f{1}{2} \max[U] \in \bR.
$$
\ethm

This result is in fact obtained by setting $\eps_m = \eps$ 
in the more general theorem appearing next.
Before that, one takes notice that the above global recovery map $\Delta^{\rm glo}$
is different, though close, to the local recovery map $\Delta^{\rm loc}$ from Theorem \ref{ThmLoc}.
Indeed, imposing $\eps_m = \eps$ for all $m \in \ival{1}{M}$, one easily arrives at the expressions
\begin{align*}
\Delta^{\rm glo}(\vy) 
&= \phantom{\f{1}{2}}\max_{m \in \ival{1}{M}} \big( y_m \big) + \f{1}{2} \max \Big[ \min_{m \in \ival{1}{M}} \big( \dist(\cdot, x^{(m)})^\alpha \big) \Big],\\
\Delta^{\rm loc}(\vy)
& = \f{1}{2} \max_{m \in \ival{1}{M}} \Big( y_m \big)
+ \f{1}{2} \max \big[ \min_{m \in \ival{1}{M}} \big( y_m + \dist(\cdot, x^{(m)})^\alpha \big) \Big],
\end{align*}
which could even be made more explicit in the univariate setting by invoking Proposition \ref{PropUnivariate}.
These expressions incidentally reveal that the global estimate is never smaller than the local estimate,
i.e., that $\Delta^{\rm glo}(\vy)  \ge \Delta^{\rm loc}(\vy)$,
as illustrated in Figure \ref{fig}(b).

\bthm
\label{ThmGloMax2}
For the recovery of the quantity of interest $f \in B(\cX) \mapsto \max[f] \in \bR$,
given the data~\eqref{DefY} and the model set \eqref{DefH},
the minimal global worst-case error is bounded below as 
$$
\inf_{\Delta: \bR^M \to \bR}
{\rm gwce}_{{\bm \eps}}(\Delta)  \ge 
\f{1}{2} \max[u_{{\bm 0}, {\bm \eps}}] + \f{1}{2} \min_{m \in \ival{1}{M}} \big( \eps_m \big),
$$
while the recovery map $\wt{\Delta} : \vy \in \bR^M \mapsto
\max_{m \in \ival{1}{M}} \big( y_m \big) + \df{1}{2} \max[u_{{\bm 0}, {\bm \eps}}] - \df{1}{2} \max_{m \in \ival{1}{M}} \big( \eps_m \big) \in \bR$
has global worst-case error bounded above as 
$$
\phantom{\inf_{\Delta: \bR^m \to Z}}
{\rm gwce}_{{\bm \eps}}(\wt{\Delta})  \le 
\f{1}{2} \max[u_{{\bm 0}, {\bm \eps}}] + \f{1}{2} \max_{m \in \ival{1}{M}} \big( \eps_m \big).
$$
\ethm

\bpf
For the lower bound, 
one points out that $\ell_{{\bm 0},{\bm \eps}} \in \cH_\alpha$ and $u_{{\bm 0},{\bm \eps}} \in \cH_\alpha$
based on an argument similar to one in the proof of Theorem \ref{ThmLoc}.
One also sees from their defining expressions that $-\eps_m \le \ell_{{\bm 0},{\bm \eps}}(x^{(m)}) \le 0 \le u_{{\bm 0},{\bm \eps}}(x^{(m)}) \le \eps_m$ for all $m \in \ival{1}{M}$,
and hence $\ell_{{\bm 0},{\bm \eps}}(x^{(m)}) + e_m = 0$ for some $\ve \in \cE_{\bm \eps}$ and $u_{{\bm 0},{\bm \eps}}(x^{(m)}) + e_m = 0$ for some $\ve \in \cE_{\bm \eps}$.
Therefore, one derives that, for any $\Delta: \bR^m \to \bR$,
\begin{align*}
{\rm gwce}_{{\bm \eps}}(\Delta)  
\ge
\begin{cases}
|\max[\ell_{{\bm 0},{\bm \eps}}] \, - \Delta(0)|\\
|\max[u_{{\bm 0},{\bm \eps}}] - \Delta(0)|
\end{cases}
& \ge \f{1}{2} \Big( |\max[u_{{\bm 0},{\bm \eps}}] - \Delta(0)|
+ |\max[\ell_{{\bm 0},{\bm \eps}}] - \Delta(0)| \Big) \\
& \ge \f{1}{2} \big( \max[u_{{\bm 0},{\bm \eps}}] - \max[\ell_{{\bm 0},{\bm \eps}}] \big).
\end{align*}
The result follows from the observation that $\max[\ell_{{\bm 0},{\bm \eps}}] = \max_{m  } \big( -\eps_m \big)
= - \min_{m  } \big( \eps_m \big)$.

For the upper bound on ${\rm gwce}_{{\bm \eps}}(\wt{\Delta})$,
one needs to establish that,
given $f \in \cH_\alpha$, $\ve \in \cE_{\bm \eps}$, and writing $\vy = {\bm \La}_{\vx} f + \ve$ for short,
$$
-\f{1}{2} \max[u_{{\bm 0}, {\bm \eps}}] - \f{1}{2} \max_{m \in \ival{1}{M}} \big( \eps_m \big)
\le
\max[f] - \wt{\Delta}(\vy)
\le
\f{1}{2} \max[u_{{\bm 0}, {\bm \eps}}] + \f{1}{2} \max_{m \in \ival{1}{M}} \big( \eps_m \big).
$$
Taking the expression of $\wt{\Delta}(\vy)$ into account in these prospective inequalities,
one notices that the leftmost one reduces to $\max[f] \ge \max_{m } \big( y_m \big) - \max_{m } \big( \eps_m \big)$,
while the rightmost one reduces to
$\max[f] \le \max_{m } \big( y_m \big) + \max[u_{{\bm 0},{\bm \eps}}]$.
The  facts that  $\max[f] \ge \max_{m } \big( f(x^{(m)}) \big)
= \max_{m } \big( y_m - e_m \big)$
and $f(x) \le u_{\vy,{\bm \eps}}(x) = \min_{m } \big( y_m + \eps_m + \dist(x,x^{(m)})^\alpha \big)
\le \max_{m } \big( y_m \big) + u_{{\bm 0},{\bm \eps}}(x)$ for all $x \in \cX$
yield the former and the latter, respectively.
\epf

\brk
Theorem \ref{ThmGloMax} could be extended in a further direction, incorporating jittered observations, so as to receive data of the 
form $y_m = f(x^{(m)}+\xi^{(m)}) + e_m$ with $|e_m| \le \eps$
and $\|\xi^{(m)} \| \le \delta$.
It is hereby assumed that this norm is the one generating the metric $\dist$ and the domain $\cX \inc \bR^d$ as its unit ball.
Defining a new notion of global error including the worst case over $\xi^{(1)},\ldots,\xi^{(M)} \in \cX$,
one would obtain ${\rm gwce}_{\eps,\delta}(\Delta) \ge \max[U]/2 + \eps$ for any $\Delta: \bR^M \to \bR$,
simply because this new notion ${\rm gwce}_{\eps,\delta}$ of global worst-case error is larger than the previous notion ${\rm gwce}_{\eps}$. 
Moreover, for the recovery map $\Delta^{\rm glo}$ from Theorem \ref{ThmGloMax},
one can obtain ${\rm gwce}_{\eps,\delta}(\Delta^{\rm glo}) \le \max[U]/2 + \eps + \delta$
by making the necessary adjustments in the second part of the proof of Theorem \ref{ThmGloMax2}. 
%(of course, one could also allow parameters $\eps_m$ and $\delta_m$ that depend on $m$).
Note that it makes sense to assume that the balls $B(x^{(m)},\delta)$ do not intersect,
otherwise a jittered observation $y_m$ could also come from another datasite $x^{(n)}$ with $n \not= m$.
Consequently, if $\mu$ is the volumetric measure, then
$$
M \delta^d \mu(B(0,1)) 
= \sum_{m=1}^M \mu(B(x^{(m)},\delta))
= \mu \Big( \bigcup_{m=1}^M B(x^{(m)},\delta) \Big) 
\le \mu (B(0,1+\delta)) \le (1+\delta)^d \mu(B(0,1)), 
$$
yielding $M^{1/d} \le 1 + 1/\delta \le 2/\delta$,
i.e., $\delta \le 2/M^{1/d}$.
Since the next section will establish that
$\max[U] \ge 1/(2M^{\alpha/d}) \ge 1/(2M^{1/d})$,
one derives that $\delta \le 4 \max[U]$,
which in turn implies that ${\rm gwce}_{\eps,\delta}(\Delta^{\rm glo}) \le 
\max[U]/2 + \eps + 4 \max[U]
\le 9 \, (\max[U]/2 + \eps) \le
9 \, {\rm gwce}_{\eps,\delta}(\Delta)$ for any $\Delta: \bR^M \to \bR$.
In other words, one can assert that $\Delta^{\rm glo}$ is near optimal with a factor $9$.
\erk

\section{Estimations of the Minimal Global Worst-Case Errors}
\label{SecExp}

Throughout this section,
given a norm $\|\cdot\|$ on $\bR^d$,
the domain $\cX \inc \bR^d$ is chosen to be its unit ball $B(0,1)$ equipped with the Lebesgue measure $\mu$ normalized so that $\mu(\cX) = 1$.
The metric involved in the H\"older model set $\cH_\alpha$ defined in \eqref{DefH} is also chosen to be the one subordinated to this norm,
i.e., $\dist(x,x') = \|x-x'\|$.
As shown in the previous section, the minimal global worst-case errors are dictated by the function $u_{{\bm 0}, {\bm \eps}}$,
e.g. via its $L_p$-norm for the recovery of $\Gamma(f)=f$,
$\Gamma$ being the identity from $B(\cX)$ to $L_p(\cX)$ (see Theorem \ref{ThmLinearQ}),
or via its $L_\infty$-norm for the recovery of $\max[f]$ (see Theorem \ref{ThmGloMax2}). 
The goal here is to evaluate these quantities optimized over the selection of datasites $x^{(1)},\ldots,x^{(M)}$.
The first result applies to a general $\cX = B(0,1)$ and gives a precise two-sided estimate valid for any $p \in [0,\infty]$.
It is followed by an exact expression applying to the case of the $\ell_\infty$-norm on $\bR^d$, so that $\cX = [-1,1]^d$.

\bthm
\label{ThmGenX}
For $\cX$ being the unit ball of $\bR^d$ relatively to a norm $\|\cdot\|$ and for $p \in [1,\infty]$, recalling that $u_{{\bm 0},{\bm \eps}}(x) = \min_{m } \big( \eps_m + \|x-x^{(m)}\|^\alpha  \big)$,
one has
$$
\min_{m \in \ival{1}{M}} \big( \eps_m \big) + \df{c}{M^{\alpha/d}} 
\le
\inf_{x^{(1)},\ldots,x^{(M)} \in \cX} \|u_{{\bm 0},{\bm \eps}}\|_{L_p(\cX)}
\le 
\max_{m \in \ival{1}{M}} \big( \eps_m \big) + \df{C}{M^{\alpha/d}} ,
$$
where the absolute constants $c,C >0$ can be taken as  $c = 1/2$ and $C = 3$.
\ethm

\bpf
For the lower bound,
given any $x^{(1)},\ldots,x^{(M)} \in \cX$,
one uses  $\|u_{{\bm 0},{\bm \eps}}\|_{L_p(\cX)} \ge \|u_{{\bm 0},{\bm \eps}}\|_{L_1(\cX)}$ and 
$u_{{\bm 0},{\bm \eps}}(x) \ge \min_{m } (\eps_m) + U(x) $
to deduce that $\|u_{{\bm 0},{\bm \eps}}\|_{L_p(\cX)} \ge \min_{m } (\eps_m) + \|U\|_{L_1(\cX)}$.
Then, in order to bound $\|U\|_{L_1(\cX)}$ from below,
one adapts the proof of \cite{Suk}, 
where Sukharev was interested in the optimal recovery of the integral of a Lipschitz function on $[0,1]^d$, but otherwise the argument is essentially the same.
Namely, for $\tau > 0$ to be chosen later, one writes
\begin{align*} 
\|U\|_{L_1(\cX)}
& = \int_{\cX} U(x) d\mu(x) 
= \int_0^\infty \mu\{ x \in \cX: U(x) >t \} dt
% \ge \int_0^t  \mu\{ x \in \cX: U(x) >t^{1/p} \} dt
\ge \int_0^\tau  \big[ 1 -  \mu\{ x \in \cX: U(x) \le t \} \big] dt\\
& = \int_0^\tau  \big[ 1 -  \mu\{ x \in \cX: \|x - x^{(m)}\|^{\alpha} \le t
\mbox{ for some } m \in \ival{1}{M} \} \big] dt\\
& \ge \int_0^\tau \Big[ 1 - \sum_{m=1}^M \mu\{ x \in \cX: \|x - x^{(m)} \|  \le t^{1/\alpha}  \} \Big] dt 
\ge \tau - \sum_{m=1}^M \int_0^\tau (t^{1/\alpha})^d  dt\\
& = \tau - \f{M}{d/\alpha + 1} \tau^{d/\alpha+1}.
\end{align*}
This lower bound is optimized when $1-M \tau^{d/\alpha}=0$, i.e., $\tau = M^{-\alpha/d}$, 
yielding
$$
\|U\|_{L_1(\cX)} \ge M^{-\alpha/d} - \f{M}{d/\alpha + 1} \big(M^{-\alpha/d}\big)^{d/\alpha + 1}
= \Big( 1 - \f{1}{d/\alpha+1} \Big) M^{-\alpha/d}
= \f{d}{d+\alpha} M^{-\alpha/d}.
$$
The announced inequality now simply follows from $d/(d+\alpha) \ge 1/2$.

For the upper bound, one uses $\|u_{{\bm 0},{\bm \eps}}\|_{L_p(\cX)} \le \|u_{{\bm 0},{\bm \eps}}\|_{L_\infty(\cX)}$ and 
$u_{{\bm 0},{\bm \eps}}(x) \le \max_{m } ( \eps_m ) + U(x) $
to deduce that $\|u_{{\bm 0},{\bm \eps}}\|_{L_p(\cX)} \le  \max_{m } ( \eps_m ) + \|U\|_{L_\infty(\cX)}$.
Then, in order to bound $\|U\|_{L_\infty(\cX)}$ from above,
one relies on a folklore result
%\footnote{The result can be found in \cite{Pis}, but it was already folklore then, dixit the author. The argument consists in considering a maximal $\theta$-separated sets, which is automatically a $\theta$-net, and bound its cardinality by a volumetric argument akin to the one presented earlier.} 
about covering numbers.
Namely, 
one can find a set $\{x^{(1)},\ldots,x^{(N)} \} \inc \cX$
of $N \le (1 + 2/\theta)^d$ points such that,
for any $x \in \cX$,
there exists $n \in \ival{1}{N}$ such that $\|x-x^{(n)}\| \le \theta$.
Selecting $\theta = 3/M^{1/d}$, one has $N \le (3/\theta)^d = M$,
so one can complete the above set to form a $\theta$-net $\{ x^{(1)}, \ldots, x^{(M)} \} \inc \cX$ consisting of $M$ points.
Using these datasites, it is now clear that $\|U\|_{L_\infty(\cX)} = \max_{x} \min_{m} \| x - x^{(m)}\|^\alpha \le (3/M^{1/d})^\alpha \le 3 / M^{\alpha/d}$.
\epf

The result of Theorem \ref{ThmGenX} can be refined for the 
 standard domain $\cX = [-1,1]^d$,
as uncovered below.
\rev{The case $p=\infty$, not included but obtainable as a limiting case, would in fact be easier to establish.}

\bthm
\label{ThmXCube}
If $\bR^d$ is equipped with the $\ell_\infty$-norm,
so that $\cX = [-1,1]^d$
and $\dist(x,x') = \|x-x'\|_{\ell_\infty^d}$,
and if all the $\eps_m$'s are equal to a common $\eps$,
then, for any $p \in [1,\infty)$,
$$
\inf_{x^{(1)},\ldots,x^{(M)} \in \cX} \|u_{{\bm 0},{\bm \eps}}\|_{L_p(\cX)}^p
\ge \left( \eps + \f{1}{M^{\alpha/d}} \right)^p - M \int_{\eps^p}^{(\eps + 1/M^{\alpha/d})^p} (t^{1/p} - \eps)^{d/\alpha} dt,
$$
with equality being achieved when $M^{1/d}$ is an integer.
Furthermore, if $\eps = 0$, then
$$
\inf_{x^{(1)},\ldots,x^{(M)} \in \cX}
\|u_{{\bm 0},{\bm \eps}}\|_{L_p(\cX)} \ge \left( \f{d}{d+\alpha p} \right)^{1/p} \f{1}{M^{\alpha/d}}
\qquad \mbox{with equality when } M^{1/d} \mbox{ is an integer}.
$$
\ethm

\bpf
For arbitrary points $x^{(1)},\ldots,x^{(M)} \in \cX$,
let $V_1,\ldots,V_M$ denote the (genuine) Voronoi cells,
i.e., $V_m = \{ x \in \cX: \|x-x^{(m)}\|_{\ell_\infty^d} \le \|x-x^{(n)}\|_{\ell_\infty^d} \mbox{ for all } n \in \ival{1}{M} \}$.
One starts by writing 
\begin{align}
\nonumber
\|u_{{\bm 0},{\bm \eps}}\|_{L_p(\cX)}^p 
& = \int_\cX \min_{m \in \ival{1}{M}} \big( \eps + \|x-x^{(m)}\|_{\ell_\infty^d}^\alpha \big)^p d\mu(x)
= \sum_{m=1}^M \int_{V_m}  \big( \eps + \|x-x^{(m)}\|_{\ell_\infty^d}^\alpha  \big)^p d\mu(x) \\
\nonumber
& = \sum_{m=1}^M \int_0^\infty \mu \{ x \in V_m: \big( \eps + \|x-x^{(m)}\|_{\ell_\infty^d}^\alpha  \big)^p >t \} dt\\
\nonumber
& = \sum_{m=1}^M \int_0^\infty \mu \{ x \in V_m:  \|x-x^{(m)}\|_{\ell_\infty^d} > (t^{1/p} -\eps)^{1/\alpha} \} dt\\
\label{UseVoronoi}
& = \sum_{m=1}^M \int_0^\infty \mu( V_m \setm B(x^{(m)}, (t^{1/p} - \eps) ^{1/\alpha}) ) dt.
\end{align}
For  $\tau>0$ to be chosen later, one continues with
\begin{align*}
\|u_{{\bm 0},{\bm \eps}}\|_{L_p(\cX)}^p 
& \ge \sum_{m=1}^M \int_0^\tau \mu( V_m \setm B(x^{(m)}, (t^{1/p} - \eps)^{1/\alpha} ) ) dt\\
& \ge \sum_{m=1}^M \int_0^\tau \big( \mu( V_m) - \mu(B(x^{(m)}, (t^{1/p} - \eps)^{1/\alpha} ) ) \big) dt\\
& = \tau  - M \int_{\eps^p}^{\tau} (t^{1/p} - \eps)^{d/\alpha} dt.
\end{align*}
Selecting the optimal $\tau = ( \eps + M^{-\alpha/d} )^p$ gives the announced lower bound.

Assuming now that $M = N^d$ for some integer $N \in \bN$,
one aims at showing that the lower bound can be achieved.
To this end, one partitions $\cX = [-1,1]^d$ into a regular grid made of $N^d$ hypercubes of side $2/N$
and one picks $x^{(1)},\ldots,x^{(M)}$ as their centers.
Thus, the hypercubes are the Voronoi cells associated with $x^{(1)},\ldots,x^{(M)}$,
i.e., $B(x^{(m)},1/N) = V_m$.
Coming back to \eqref{UseVoronoi}
while remarking that 
$$
\mu( V_m \setm B(x^{(m)}, (t^{1/p} - \eps) ^{1/\alpha}) )
= \begin{cases}
N^{-d}, & \mbox{if } t \le \eps^p,\\
N^{-d} - (t^{1/p} - \eps) ^{d/\alpha}, & \mbox{if } \eps^p < t \le (\eps + N^{-\alpha})^p,\\
0, & \mbox{if } t> (\eps + N^{-\alpha})^p,\\
\end{cases}
$$
one obtains, for this choice of $x^{(1)},\ldots,x^{(M)}$,
\begin{align*}
\|u_{{\bm 0},{\bm \eps}}\|_{L_p(\cX)}^p
& = \sum_{m=1}^M \Big(
\int_0^{\eps^p} \f{1}{N^d} dt + \int_{\eps^p}^{(\eps+N^{-\alpha})^p} \Big( \f{1}{N^d} - (t^{1/p} - \eps) ^{d/\alpha} \Big) dt
\Big)\\
& = (\eps+N^{-\alpha})^p
- M \int_{\eps^p}^{(\eps+N^{-\alpha})^p}  (t^{1/p} - \eps) ^{d/\alpha}  dt.
\end{align*}
Substituting $N = M^{1/d}$ shows that the latter indeed equals the previous lower bound.  
\epf

With $p=1$ and $\alpha = 1$,
the infimum from Theorem \ref{ThmXCube} becomes $\eps + (d/(d+1)) M^{-1/d} $.
Similarly,
Theorem \ref{ThmGenX} also demonstrated that in general the minimal global-worst case error over all possible recovery procedures and datasites is composed of two terms:
$\eps$, the observation accuracy,
and $M^{-1/d}$, where $M$ is the number of observations.
The term $M^{-1/d}$ is bad news:
even for perfectly accurate observations,
an error of order $\eta<1$ can only be achieved with a number~$M$ of observations of order $\eta^{-d}$---this is the curse of dimensionality, as formalized in Information-Based Complexity.
But the high dimensionality seems to come with a silver lining, too:
\rev{in a scenario when the number $M$ of observations is fixed, independent of $d$,
and when the observation accuracy $\eps$ can be selected,
the smallest error achieved at $\eps = 0$
is of the same order $M^{-1/d}$ as if one chose $\eps = M^{-1/d}$, which is less stringent as the dimension $d$ grows.}

%at a fixed $M$,
%since the error contains the term $M^{-1/d}$ anyway,
%it is pointless to observe with accuracy~$\eps$ better than the  threshold~$M^{-1/d}$,
%which is less stringent as the dimension $d$ grows.

\end{document}